\documentclass[preprint,subeqn,12pt]{elsarticle}
\usepackage{latexsym}
\usepackage{amsfonts}
\usepackage{amsmath}
\usepackage{epsf}
\usepackage{graphicx}
\usepackage{amssymb}
\usepackage{pifont}
\usepackage{fleqn}
\usepackage{txfonts}
\usepackage{hyperref}
\usepackage{lineno}

\newdefinition{remark}{Remark}
\newdefinition{prop}{Proposition}

\begin{document}

\begin{frontmatter}

\title{Linear dynamic model of production-inventory with debt repayment: optimal management strategies}
\author[Ekaterina]{Ekaterina Tuchnolobova}
\ead{tuchnokaterina@mail.ru}
\author[Victor]{Victor Terletskiy}
\ead{terletsky@math.isu.ru}
\address[Ekaterina,Victor]{Institute for Mathematics, Economics and Informatics, Irkutsk State University, Bulvar Gagarina 20, Irkutsk 664003, Russia}
\author[Olga]{Olga Vasilieva\corref{cor1}}
\ead{olgavas@univalle.edu.co}
\cortext[cor1]{Corresponding author}
\address[Olga]{Department of Mathematics, Universidad del Valle, Calle 13 No. 100-00, Cali, Colombia }

\begin{abstract}
In this paper, we present a simple microeconomic model with linear continuous-time dynamics that describes a production-inventory system with debt repayment. This model is formulated in terms of optimal control and its exact solutions are derived by prudent application of the maximum principle under different sets of initial conditions (scenarios). For a potentially profitable small firm, we also propose some alternative short-term control strategies resulting in a positive final profit and prove their optimality. Practical implementation of such strategies is also discussed.
\end{abstract}

\begin{keyword}
production-inventory system \sep debt repayment \sep linear dynamics \sep optimal control \sep profitability condition

\MSC 91B55 \sep 49N90 \sep 91B38
\end{keyword}

\end{frontmatter}

\section{Introduction}
\label{sec:intro}

In macro- and microeconomics, the continuous-time models involving ordinary differential equations naturally serve as a basis for understanding the behavior of economic systems where the dynamic aspects play an important role.

Mathematical models of microeconomic can help to explain macroeconomic phenomena and to improve the management of a particular production unit (plant, firm, family business, etc.). Among them, the major focus has been placed on \emph{production-inventory systems} consisting of a manufacturing plant and a warehouse to store the finished goods which are produced but not immediately sold. Usually, the production rate is treated as a control variable while the purpose of control consists in meeting the existent level of product demand at the market or in maximizing the net profit of the production unit.

Traditional approach to solving the production-inventory problem in terms of optimal control is lucidly described in \cite{SethiThompson2000} while a state-of-art review \cite{OrtegaLin2004} identifies the major research efforts for applying control theoretic methods to production–inventory systems. During the last decade, various scholars
have made essential contributions to this field. Among them it is worth to mention the control theory applications to the models  with stock-dependent demand rate \cite{TengChang2005} and with inventory-level-dependent demand rate \cite{KhmelnitskyGerchak2002,Urban2005}. Other stream of research has been focused on optimization of inventory systems with product deterioration (see, e.g., \cite{Al-khedhairiTadj2006,BatenKamil2011,BakkerRiezebosTeuner2012}), systems with back-orders and lost sales \cite{ChangLo2009,BenjaafarElHafsiHuang2010} or without them \cite{AdidaPerakis2007}.

However, all previously mentioned works do not consider explicitly the dynamics of the firm's debts acquired during the production period or prior to its commencement (overdue payables).

On the other hand, V. Tokarev \cite{Tokarev2001,Tokarev2002} had proposed a microeconomic model of short-term crediting and debt repayment for a small firm where the rate of debt repayment has been treated as a control variable. This model does not consider inventory stock variation and mainly concentrates on the dynamic of production funds (firm's basic assets) and current debts. The purpose of control consists in maximizing the final value of the firm's production funds and the optimal control has a \emph{bang-bang} structure naturally dependent on the loan's interest rate. Namely, for low-interest banking rates it is more profitable to invest all available cash into production at the beginning of the period (in order to generate more profit) and then to repay the debts by the end of the period. For high-interest banking rates a reverse strategy is optimal: first comes the debt repayment and then the investment, since in this case the debts grow faster than the maximal possible profit obtained from production. The same model was further developed in \cite{GrigorievaKhailov2005} under additional condition that all pending debts must be paid off by the end of the period.

In this paper, we propose a simple microeconomic model with linear continuous-time dynamics that explicitly includes the variation of current debt of a firm. Our model combines the traditional features of production-inventory systems with Tokarev's approach and has three constrained state variables together with three bounded control variables. The model is premeditated for dealing with \emph{short-term} planning of a \emph{profitable} firm (such as small or family business) in \emph{economically stable} market environment and may provide some guidelines for periodic short-term planning (weekly, monthly, etc.).

Linear structure of the model allows to apply the Pontryagin's maximum principle as a \emph{necessary and sufficient} condition of optimality. However, the application of the maximum principle for problems with mixed-type constraints can be rather challenging (see, e.g. \cite{Maurer1977,HartlSethiVickson1995,MaurerPesch2008}). Sometimes, there is no other way than to ``guess'' a possibly optimal control strategy and then to use the maximum principle for proving its optimality. For such ``guesses'', we have just used a ``common-sense'' control strategies derived from the economic contexts of the problem. We have deliberately kept this model simple in order to obtain its analytical solutions under different scenarios and to propose some viable alternatives for better management in economically stable environment.

The paper is organized as follows. Section \ref{sec:model} provides the model description and formulates the maximum principle for optimal control problems with mixed-type constraints. Section \ref{sec:scenarios} presents three basic scenarios of the firm management for different sets of initial conditions. In Section \ref{sec:alter} we discuss some alternative strategies of the firm's management and justify their optimality as well as practical implementation within the frameworks of the model. Section \ref{sec:conclusion} contains the conclusion and briefly indicates the perspectives of future research.

\section{Model description and preliminaries}
\label{sec:model}

In economically stable market environment, the accumulation of arrears by a small firm (family business, for example) clearly indicates the inefficiency of its short-term planning. Therefore, financial manager or self-employed entrepreneur should try to avoid the accumulation of unpaid debts while seeking to maximize the firms' cumulative profit.

To examine this situation we propose to use a simple dynamic model of a firm producing a single good which operates under market stability\footnote{Under market stability, the inflation index for short-term periods is considered to be depreciable. Therefore, we omit the intertemporal discount factor in the model setting. } and has a viable (deterministic) estimation for its product demand at the market.  Such model can be posed in terms of optimal control as follows:

\begin{equation}
\label{objective}
\mathcal{J} ( u,v,w ) = N(T) - D(T) \to \max
\end{equation}
subject to
\begin{subequations}
\label{system}
    \begin{equation}
    \label{system-N}
    \dot{N}(t) = p w(t) - v(t) - Z(u), \quad N(0)=N_0,
    \end{equation}
\begin{equation}
\label{system-D}
\dot{D}(t) = r D(t) + A u(t) - v(t), \quad D(0)=D_0,
\end{equation}
    \begin{equation}
    \label{system-S}
    \dot{S}(t) = u(t) - w(t) - \alpha S(t), \quad S(0)=S_0,
    \end{equation}
\end{subequations}
under state constraints
\begin{equation}
\label{state-constr}
N(t) \geq 0, \quad D(t) \geq 0, \quad 0 \leq S(t) \leq S_{max}, \quad t \in [0,T]
\end{equation}
and control constraints
\begin{equation}
\label{cont-constr}
0 \leq u(t) \leq u_{max}, \quad 0 \leq v(t) \leq v_{max}, \quad 0 \leq w(t) \leq w_{max}, \quad t \in [0,T].
\end{equation}
The quantities used in \eqref{objective}-\eqref{cont-constr} are defines as

\noindent
\begin{tabular}{rcl}
   $u(t)$ & \!--\! & production rate (control variable) with maximum gross output  \\
          &    & given by $u_{max}>0$ (constant); \\
  $v(t)$ & \!--\! & rate of the debt repayment (control variable) with maximum \\
         &    & repayment capacity given by $v_{max} >0$ (constant);  \\
  $w(t)$ &\!--\! & rate of sales (control variable) with maximum volume of demand \\
         &    &  given by $w_{max}>0$ (constant); \\
  $N(t)$ &\!--\! & cumulative net profit of the firm by time $t$ (state variable); \\
  $D(t)$ &\!--\! & amount of overdue payables (debts) of the firm at time $t$ (state \\
         &    & variable); \\
  $S(t)$ & \!--\! & volume of finished goods inventory at time $t$ (state variable) with \\
         &    & maximum storage capacity given by $S_{max} >0$ (constant); \\
  $Z(u)$ &\!--\! & cost function, which includes other expenditures not related to the \\
         &    & purchase of raw materials ( such as wages, social security costs, lease \\
         &    & payments, etc.). For simplicity sake, we suppose that $Z(u)=Ku + B$ \\
         &    & where $K, B >0 $ are given constants; \\
  $p>0$ & \!--\! & output (retail) price index (constant); \\
  $r>0$ & \!--\! & rate of overdue debt accumulation (constant); \\
  $\!\! A>0$ & \!--\! & average cost for the purchase of raw materials (constant); \\
  $\alpha>0$ & \!--\! & outflow rate of finished goods inventory (constant) that includes sales \\
        &     & and stock loss during storage.
\end{tabular}
\noindent
Here the maximum rate of the debt repayment $v_{max}$ is supposed to be sufficiently large while the maximum volume of demand $w_{max}$ does not naturally exceed the maximum level of production capacity $u_{max}$, that is, $w_{max} \leq u_{max} $. It should be noted that $Z(u)=Ku + B$ is chosen as \emph{affine} linear function in order to emphasize that even for $u=0$, these indirect costs will not be zero.

Criterion \eqref{objective} expresses the maximization of the cumulative net profit $N(t)$ and minimization of overdue debts $D(t)$ by final time $T$. Here we do not impose the condition of \emph{total} debt repayment by the final time (that is, $D(T)=0$ as in \cite{GrigorievaKhailov2005}) and merely try to minimize the debt that may remain positive.

Equation \eqref{system-N} describes the net profit accumulation: total revenue $p w(t)$ minus the debt repayment $v(t)$ and other costs $Z(u)$. Equation \eqref{system-D} provides the dynamics of accumulated arrears: the debt increases due to the interest rate $r$ and purchase of raw materials $Au(t)$ and decreases with repayments $v(t)$. Equation \eqref{system-S} describes changes in finished goods inventory: the stock increases with production and decreases with sales and losses during storage.

By introducing formal notation $\mathbf{X}(t) = \left( N(t),D(t),S(t) \right)^{\prime}, \mathbf{U}(t)=\left( u(t), v(t), w(t) \right)^{\prime}$ and $\mathbf{\Psi} = \left( \psi_1, \psi_2, \psi_3 \right)^{\prime}$ (vector of adjoint variables associated with state variables $N(t),D(t),$ and $S(t)$, respectively), the Hamiltonian of the optimal control problem \eqref{objective}--\eqref{cont-constr} is linear in $\mathbf{U}$ and separable with respect to $ \left( u(t), v(t), w(t) \right)$:
\[ H \left( \mathbf{\Psi},\mathbf{X},\mathbf{U} \right) = \psi_1 \left[ p w - v - K u - B \right] + \psi_2 \left[ r D + A u - v \right] + \psi_3 \left[ u - w - \alpha S \right]
\]
\[ = \left[ -K \psi_1 + A \psi_2 + \psi_3 \right] u - \left[ \psi_1 + \psi_2 \right] v + \left[ p \psi_1 - \psi_3 \right] w -B \psi_1 + r D \psi_2 - \alpha S \psi_3, \]
where the factors to the control components are the switching functions:
\begin{equation}
\label{switch}
\theta_u = -K \psi_1 + A \psi_2 + \psi_3, \quad \theta_v = - \psi_1 - \psi_2, \quad \theta_w = p \psi_1 - \psi_3.
\end{equation}
Note that the switching vector $ \mathbf{\Theta}= \left( \theta_u,\theta_v,\theta_w \right)$ only depends on the adjoint vector $\mathbf{\Psi} $ (whose components are also called ``shadow prices" in economics). By assigning a vector of Lagrange multipliers $\mathbf{\Lambda} = \left( \lambda_1, \lambda_2,\lambda_3,\lambda_4\right)^{\prime}$ to  four state constraints \eqref{state-constr}, we can define the Lagrangian for our problem as
\begin{equation}
\label{lagrangian}
\mathcal{L} \left( \mathbf{\Psi},\mathbf{X},\mathbf{U}, \mathbf{\Lambda} \right) = H \left( \mathbf{\Psi},\mathbf{X},\mathbf{U} \right) + \lambda_1 N + \lambda_2 D + \lambda_3 S + \lambda_4 \left[ S - S_{max} \right] \nonumber
\end{equation}
For our linear control problem \eqref{objective}--\eqref{cont-constr}, the maximum principle serves as \emph{necessary and sufficient condition} for optimality and can be rigorously justified using \emph{direct adjoining approach} described, e.g., in \cite{Maurer1977,HartlSethiVickson1995,MaurerPesch2008}. Therefore, a piecewise continuous (or \emph{bang-bang}) control $ \mathbf{U}^{*}(t)=\left( u^{*}, v^{*}(t), w^{*}(t) \right)^{\prime} $ defined by switching functions \eqref{switch} that maximizes the Hamiltonian almost in all points of $[0,T]$ is optimal in \eqref{objective}--\eqref{cont-constr} if and only if there exist an absolutely continuous costate trajectory  $\mathbf{\Psi}: [0,T] \to \mathbb{R}^3$ of the adjoint system
\begin{eqnarray}
\dot{\psi}_1(t) = - \frac{\partial \mathcal{L}}{\partial N} = - \lambda_1(t), & & \psi_1(T) = \mu_1 + 1, \nonumber \\
\label{adsystem}
\dot{\psi}_2(t) = - \frac{\partial \mathcal{L}}{\partial D} = - r \psi_2(t) - \lambda_2(t), & & \psi_2(T) = \mu_2 - 1, \\
\dot{\psi}_3(t) = - \frac{\partial \mathcal{L}}{\partial S} = \alpha \psi_3(t) - \lambda_3(t) + \lambda_4(t), & & \psi_3(T) = \mu_3 - \mu_4, \nonumber
\end{eqnarray}
as well as a piecewise continuous vector function of multipliers $\mathbf{\Lambda}: [0,T] \to \mathbb{R}^4$, and a nonzero vector $\mathbf{\mu}=\left( \mu_1, \mu_2, \mu_3, \mu_4 \right)^{\prime}$ such that the following conditions of \emph{complementary slackness} are satisfied:
\begin{equation}
\label{slackness}
\begin{array}{rr}
\lambda_1(t) N(t) = 0 & \mu_1 N(T) = 0 \\
\lambda_2(t) D(t) = 0 & \mu_2 D(T) = 0 \\
\lambda_3(t) S(t) = 0 & \mu_3 S(T) = 0 \\
\lambda_4(t) \left[ S(t) - S_{max} \right] = 0 & \mu_4 \left[ S(T) - S_{max} \right]  = 0 \end{array} \quad
\begin{array}{c}
\lambda_i(t) \geq 0 \\ \mu_i \geq 0, \\ i=1,2,3,4.
\end{array}
\end{equation}
This classical result will be very essential for further design of optimal control strategies $\mathbf{U}^{*}= \left( u^{*},v^{*},w^{*} \right)^{\prime}$ under which the objective functional \eqref{objective} attains its maximum value.

\section{Optimal functioning of a profitable firm: case studies.}
\label{sec:scenarios}

Optimal functioning of a firm significantly depends on whether an external demand for its products ensures a positive profit. In mathematical formalization (see, e.g. \cite{Kessides1990}), a firm is profitable if the following inequality holds:
\begin{equation}
\label{profit-con}
p w(t) > Z(u) + A u(t),
\end{equation}
that is, if its sales profit fully covers all underlying expenses (such as purchase of raw materials, equipment, and other indirect costs). Otherwise, the firm is unprofitable.

Generally speaking, optimal management strategies will also depend on the initial values of state variables, such as presence or absence of initial profits, debts, and finished goods inventory. Therefore, it is interesting from the economic point of view to revise three basic scenarios and obtain optimal strategies satisfying the maximum principle. In the subsequent case-studies the condition \eqref{profit-con} will be in force.

\subsection{Scenario 1: absence of initial debt, presence of stock at $t=0$}
\label{subsec31}

Suppose that at initial time $t=0$ the firm has no arrears and possesses  some stock of finished goods  and cash resources, that is,
\[ N(0)=N_0 > 0, \quad D(0)=D_0 =0, \quad S(0)=S_0 > 0.  \]
In this case, it will be profitable for the firm to start the production at the moment $t_S \in (0,T)$ by which the whole existing stock of finished goods is sold, thus generating no new debt before $t_S$. The production is then started at the moment $t_S$ with maximum volume of demand $w_{max}$ in order to avoid overproduction. Then the production costs of the firm at any given time $t \in [t_S, T]$ will be $A w_{max} \leq v_{max}$ (since $v_{max}$ is rather large).  Therefore, we anticipate that optimal control will be of the form:
\begin{equation}
\label{oc-scenario1}
\hspace{-5mm} u^{*} = \left\{ \begin{array}{rl} 0, & t \in [0, t_S) \\
w_{max}, & t \in [t_S,T] \end{array} \right., \quad  v^{*} = \left\{ \begin{array}{rl} 0, & t \in [0, t_S) \\
Aw_{max}, & t \in [t_S,T] \end{array} \right., \quad w^{*} = w_{max}.
\end{equation}
Under this strategy, no arrears will arise (that is, $D(t) \equiv 0, t \in [0,T]$), because the firm is able to pay its debts on time, while interest payments on the debt are economically unprofitable. In order to prove the optimality of \eqref{oc-scenario1} and find a point $t_S \in (0,T)$, we must determine Lagrange multipliers $\lambda_i(t) \geq 0$ and $\mu_i \geq 0, i=1,2,3,4$ that satisfy the complementary slackness conditions \eqref{slackness}:
\[ N(t) > 0 \quad \mbox{implies} \quad \lambda_1(t) = 0, \; \mu_1 = 0; \]
\[ D(t) = 0 \quad \mbox{implies} \quad \lambda_2(t) \geq 0, \; \mu_2 \geq 0; \]
\[ S(t) \left\{ \begin{array}{rl} > 0, & t \in [0, t_S) \\ = 0, & t \in [t_S, T] \end{array} \right.  \quad \mbox{implies} \quad \lambda_3(t) \left\{ \begin{array}{rl} = 0, & t \in [0, t_S) \\ \geq 0, & t \in [t_S, T] \end{array} \right., \; \mu_3 \geq 0; \]
\[ S(t) < S_{max} \quad \mbox{implies} \quad \lambda_4 (t) = 0, \; \mu_4 = 0. \]
The underlying adjoint system \eqref{adsystem} can be written as
\begin{equation}
\label{adsys-sc1}
\begin{array}{rclcl}
\psi_1(t) & \equiv & 1, & & \\
\dot{\psi}_2(t) & = & - r \psi_2(t) - \lambda_2(t), & & \psi_2(T) = \mu_2 - 1, \\
\dot{\psi}_3(t) & = & \left\{ \begin{array}{ll} \alpha \psi_3(t), & t \in [0, t_S) \\ \alpha \psi_3(t) - \lambda_3(t), & t \in [t_S,T] \end{array} \right., & & \psi_3(T) = \mu_3, \end{array}
\end{equation}
and the switching functions \eqref{switch} must satisfy the conditions:
\begin{eqnarray}
\label{switch-sc1}
\theta_u (t) = A \psi_2(t) + \psi_3(t) - K \left\{ \begin{array}{rl} < 0, & t \in [0, t_S), \\
= 0, & t \in [t_S,T], \end{array} \right. & & \\
\theta_v (t) = - \psi_2 - 1 \left\{ \begin{array}{rl} < 0, & t \in [0, t_S), \\
= 0, & t \in [t_S,T], \end{array} \right. & & \theta_w (t) = p - \psi_3(t) > 0. \nonumber
\end{eqnarray}
According to \eqref{switch-sc1}, for $t \in [t_S,T]$ it is fulfilled that
\[  A \psi_2(t) + \psi_3(t) - K = 0, \quad \psi_2(t) = -1 \]
and in view of \eqref{adsys-sc1} we have $\lambda_2(t) = r >0, \mu_2 = 0.$ On the other hand, using $\psi_1(t) \equiv 1, \psi_2(t) \equiv -1$ in \eqref{switch-sc1} it is obtained that $\psi_3(t) = A + K, t \in [t_S,T]$. By substituting this expression in the last equation of \eqref{adsys-sc1} it becomes clear that $ \lambda_3(t) = \alpha (A + K) >0$ for $ t \in [t_s,T]$ and $\mu_3 = A + K > 0$. Thus, we have found a set of multipliers $\lambda_i(t) \geq 0, \mu_i \geq 0, i = 1,2,3,4$ that satisfy the complementary slackness conditions \eqref{slackness}:
\[ \begin{array}{cccc}
\lambda_1(t) = 0, & \lambda_2(t) = r, & \lambda_3(t) = \left\{ \begin{array}{ll}  0, & t \in [0,t_S) \\ \alpha (A + K), & t \in [t_S,T] \end{array} \right., & \lambda_4(t) = 0 \\
\mu_1 = 0, & \mu_2 = 0, & \mu_3 = A + K, & \mu_4 = 0. \end{array} \]
The latter proves the optimality of control \eqref{oc-scenario1}. In order to find the switching point $t_S \in (0,T),$ let us consider the ODE system \eqref{system} under optimal control \eqref{oc-scenario1} within the interval $t \in [0, t_S)$:
\[ \begin{array}{ll}
\dot{N}(t) = p w_{max} - B, & N(0)=N_0 > 0, \\
\dot{D}(t) = rD(t), & D(0)=0, \\
\dot{S}(t) = - \alpha S(t) - w_{max}, & S(0)=S_0 >0.
\end{array} \]
whose solution is given by
\[ \hspace{-10mm} N(t) = N_0 + (p w_{max} - B)t, \quad D(t) = 0, \quad S(t) = \frac{\alpha S_0 + w_{max}}{\alpha} \exp[-\alpha t] - \frac{w_{max}}{\alpha}. \]
Apparently, the switching point $t_S$ must be a \emph{unique} root of equation $S(t)=0$ where $S(t)$ is given above, since this real function is strictly decreasing. Therefore, exactly by the moment
\begin{equation}
\label{t-s-sc1}
t_S = \frac{1}{\alpha} \ln \frac{\alpha S_0 + w_{max}}{w_{max}} > 0
\end{equation}
the firm's stock of finished goods becomes empty.
\begin{remark}
\label{rem1}
In effect, if it occurs that $t_S \geq T$ in \eqref{t-s-sc1}, then optimal control \eqref{oc-scenario1} simply becomes
\[ u^{*} = 0, \quad v^{*} = 0, \quad w^{*} = w_{max} \]
and implies no production, only sales of existent stock of finished goods until final time $T$. This situation may arise when the rate of outflow of finished goods inventory $\alpha$ is very slow while initial stock $S_0$ is replete; in other words, when  $S_0 > w_{max} (\exp[\alpha T] - 1)/\alpha.$
\end{remark}
To evaluate the objective functional \eqref{objective} (whose value is solely defined by $N(T)$ in this case) in optimal control \eqref{oc-scenario1}, we should find the solution of the corresponding ODE system \eqref{system} within the interval $t \in [t_S, T]:$
\[ \begin{array}{ll}
\dot{N}(t) = (p - A - K)w_{max} - B, & N(t_S)=N_0+ (pw_{max} - B)t_S, \\
\dot{D}(t) = rD(t), & D(t_S)=0, \\
\dot{S}(t) = - \alpha S(t), & S(t_S)= 0.
\end{array} \]
Its solution is
\[ N(t) = N_0 + (p w_{max} - B)t + (A + K)w_{max} (t_S - t), \quad D(t) = 0, \quad S(t) = 0 \]
and yields
\begin{equation}
\label{objective-sc1}
\mathcal{J} \left( u^{*},v^{*},w^{*} \right) = N_0 + (p w_{max} - B)T + (A + K)w_{max} (t_S - T).
\end{equation}

\subsection{Scenario 2: presence of initial debt and stock at $t=0$}
\label{subsec32}

Suppose that at initial time $t=0$ the firm has non-zero arrears and possesses some stock of finished goods  and cash resources, that is,
\[ N(0)=N_0 > 0, \quad D(0)=D_0 > 0, \quad S(0)=S_0 > 0.  \]
In this case, it will be profitable for the firm to start immediately the debt repayment and avoid further accumulation of arrears. The latter can be done by changing the second component of control strategy \eqref{oc-scenario1} resulting in
\begin{equation}
\label{oc-scenario2}
\hspace{-4mm} u^{*} = \left\{ \begin{array}{rl} 0, & t \in [0, t_S) \\
w_{max}, & t \in [t_S,T] \end{array} \right., \quad  v^{*} = \left\{ \begin{array}{rl} v_{max}, & t \in [0, t_D) \\
Aw_{max}, & t \in [t_D,T] \end{array} \right., \quad w^{*} = w_{max}.
\end{equation}
Here $t_S \in (0,T)$ has the same meaning as before (that is, $S(t_S)=0$) and $t_D \in (0,T)$ indicates the moment of full repayment of arrears (that is, $D(t_D)=0$).

Lagrange multipliers $\lambda_i(t) \geq 0$ and $\mu_i \geq 0, i=1,2,3,4$ must satisfy the complementary slackness conditions \eqref{slackness}:
\[ N(t) > 0 \quad \mbox{implies} \quad \lambda_1(t) = 0, \; \mu_1 = 0; \]
\[ D(t) \left\{ \begin{array}{rl} > 0, & t \in [0, t_D) \\ = 0, & t \in [t_D, T] \end{array} \right.  \quad \mbox{implies} \quad \lambda_2(t) \left\{ \begin{array}{rl} = 0, & t \in [0, t_D) \\ \geq 0, & t \in [t_D, T] \end{array} \right., \; \mu_2 \geq 0; \]
\[ S(t) \left\{ \begin{array}{rl} > 0, & t \in [0, t_S) \\ = 0, & t \in [t_S, T] \end{array} \right.  \quad \mbox{implies} \quad \lambda_3(t) \left\{ \begin{array}{rl} = 0, & t \in [0, t_S) \\ \geq 0, & t \in [t_S, T] \end{array} \right., \; \mu_3 \geq 0; \]
\[ S(t) < S_{max} \quad \mbox{implies} \quad \lambda_4 (t) = 0, \; \mu_4 = 0. \]
The underlying adjoint system \eqref{adsystem} can be written as
\[ \begin{array}{rclcl}
\psi_1(t) & \equiv & 1, & & \\
\dot{\psi}_2(t) & = & \left\{ \begin{array}{ll} - r \psi_2(t), & t \in [0, t_D) \\ -r \psi_2(t) - \lambda_2(t), & t \in [t_D,T] \end{array} \right., & & \psi_2(T) = \mu_2 - 1, \\
\dot{\psi}_3(t) & = & \left\{ \begin{array}{ll} \alpha \psi_3(t), & t \in [0, t_S) \\ \alpha \psi_3(t) - \lambda_3(t), & t \in [t_S,T] \end{array} \right., & & \psi_3(T) = \mu_3,
\end{array} \]
and the switching functions \eqref{switch} must satisfy the conditions:
\begin{eqnarray*}
\theta_u (t) = A \psi_2(t) + \psi_3(t) - K \left\{ \begin{array}{rl} < 0, & t \in [0, t_S), \\
= 0, & t \in [t_S,T], \end{array} \right. & & \\
\theta_v (t) = - \psi_2 - 1 \left\{ \begin{array}{rl} < 0, & t \in [0, t_D), \\
= 0, & t \in [t_D,T], \end{array} \right. & & \theta_w (t) = p - \psi_3(t) > 0.
\end{eqnarray*}
By employing a technique similar to the one used in the analysis of Scenario \ref{subsec31} and considering two cases ($t_S > t_D$ and $t_S < t_D$), we can find a set of multipliers  $\lambda_i(t) \geq 0$ and $\mu_i \geq  0, I = 1,2,3,4$ that satisfy the complementary slackness condition \eqref{slackness}:
\[ \hspace{-10mm} \begin{array}{ccc}
\lambda_1(t) = \lambda_4(t)= 0, & \lambda_2(t) = \left\{ \begin{array}{ll}  0, & t \in [0,t_D) \\ r, & t \in [t_D,T] \end{array} \right. \!\!\!, & \lambda_3(t) = \left\{ \begin{array}{ll}  0, & t \in [0,t_S) \\ \alpha (A + K), & t \in [t_S,T] \end{array} \right.\!\!\!,  \\
\mu_1 = \mu_4 = 0, & \mu_2 = 0, & \mu_3 = A + K. \end{array} \]
The latter proves the optimality of control \eqref{oc-scenario2}. However, it is not quite clear which switching time ($t_S$ or $t_D$) will occur first. Common sense suggests that smaller initial debt $D_0$ can be repaid faster. Therefore, if $D_0$ is relatively small, then $t_D < t_S$; otherwise, $t_D > t_S$ for relatively large $D_0$. Eventually, it may also happen that $t_D = t_S$ and there will only one switching point. The following proposition summarizes this idea and provides exact formulae for calculation of $t_D$ and $t_S$ in terms of problem entries, as well as their position with respect to each other.
\begin{prop}
\label{prop1}
The stock of finished goods becomes empty at the moment $t_S$ given by \eqref{t-s-sc1} independently of initial debt amount $D_0 >0$. The total debt repayment occurs at the moment $t_D$ that depends on $D_0$ in the following way: \\
(a) if $D_0 = v_{max} \left( 1 - \exp[- rt_S] \right)/r$ then $t_D = t_S$ and there is only one switching point; \\
(b) if $D_0 < v_{max} \left( 1 - \exp[- rt_S] \right)/r$ then $t_D < t_S$ and
\begin{equation}
\label{t-d-sc2a}
t_D = \frac{1}{r} \ln \left( \frac{v_{max}}{v_{max} - r D_0} \right);
\end{equation}
(c) if $D_0 > v_{max} \left( 1 - \exp[- rt_S] \right)/r$ then $t_D > t_S$ and
\begin{equation}
\label{t-d-sc2b}
t_D = \frac{1}{r} \ln \left( \frac{v_{max} - A w_{max}}{v_{max} - r D_0 -A w_{max} \exp[-r t_S]} \right).
\end{equation}
\end{prop}
Formal proof of this proposition can be consulted in the Appendix.
\begin{remark}
\label{rem2}
In effect, if it occurs that $t_D \geq T$, then optimal control \eqref{oc-scenario2} simply becomes
\[  u^{*} = \left\{ \begin{array}{rl} 0, & t \in [0, t_S), \\
w_{max}, & t \in [t_S,T] \end{array} \right. \quad  v^{*} = v_{max}, \quad w^{*} = w_{max}. \]
and implies constant debt repayment at maximum rate for all $t \in [0,T].$ In this case, initial debt $D_0$ must be very large: $D_0 > v_{max} \left( 1 - \exp[- rT] \right)/r$ and, therefore, current debts $D(t)$ will not be repayed by final time $T$. Additionally, Remark \ref{rem1} remains valid under this scenario as well.
\end{remark}
Otherwise, if $t_D < T$, then $D(T)=0$ regardless of position of $t_D$ with respect to $t_S$ and the value of criterion \eqref{objective} is solely defined by the cumulative net profit at final time $T$. The following proposition provide an explicit formula for evaluation of the objective functional.
\begin{prop}
\label{prop2}
For $t_D < T$ and regardless of its position with respect to $t_S$, we have
\begin{equation}
\label{objective-sc2}
\hspace{-10mm} \mathcal{J} ( u^{*},v^{*},w^{*} ) = N_0 + (A w_{max} - v_{max})t_D + K w_{max} t_S  + w_{max} (p - A -K) T - BT,
\end{equation}
where $t_S$ is given by \eqref{t-s-sc1} and $t_D$ is defined either by \eqref{t-d-sc2a} or by \eqref{t-d-sc2b} according to Proposition \ref{prop1}.
\end{prop}
The proof of this proposition is rather straightforward and its key features a given in the Appendix.

\subsection{Scenario 3: presence of initial debt and absence of initial stock at $t=0$}
\label{subsec33}

Suppose that at initial time $t=0$ the firm possess some cash resources and has no stock of finished good along with non-zero arrears, that is,
\[ N(0)=N_0 > 0, \quad D(0)=D_0 > 0, \quad S(0)= S_0 = 0.  \]
In this case, it will be profitable to start the production immediately and to pay off the existing debts straightaway with maximum rate of repayment. Therefore, we anticipate that optimal control will be of the form:
\begin{equation}
\label{oc-scenario3}
u^{*} = w_{max}, \quad  v^{*} = \left\{ \begin{array}{rl} v_{max}, & t \in [0, t_D) \\
Aw_{max}, & t \in [t_D,T] \end{array} \right., \quad w^{*} = w_{max}.
\end{equation}
Here $t_D \in (0,T)$ indicates the moment of full repayment of all debts (that is, $D(t_D)=0$).

Lagrange multipliers $\lambda_i(t) \geq 0$ and $\mu_i \geq 0, i=1,2,3,4$ must satisfy the complementary slackness conditions \eqref{slackness}:
\[ N(t) > 0 \quad \mbox{implies} \quad \lambda_1(t) = 0, \; \mu_1 = 0; \]
\[ D(t) \left\{ \begin{array}{rl} > 0, & t \in [0, t_D) \\ = 0, & t \in [t_D, T] \end{array} \right.  \quad \mbox{implies} \quad \lambda_2(t) \left\{ \begin{array}{rl} = 0, & t \in [0, t_D) \\ \geq 0, & t \in [t_D, T] \end{array} \right., \; \mu_2 \geq 0; \]
\[ S(t) = 0 \quad \mbox{implies} \quad \lambda_3(t) \geq 0, \; \mu_3 \geq 0; \]
\[ S(t) < S_{max} \quad \mbox{implies} \quad \lambda_4 (t) = 0, \; \mu_4 = 0. \]
The underlying adjoint system \eqref{adsystem} can be written as
\[ \begin{array}{rclcl}
\psi_1(t) & \equiv & 1, & & \\
\dot{\psi}_2(t) & = & \left\{ \begin{array}{ll} - r \psi_2(t), & t \in [0, t_D) \\ -r \psi_2(t) - \lambda_2(t), & t \in [t_D,T] \end{array} \right., & & \psi_2(T) = \mu_2 - 1, \\
\dot{\psi}_3(t) & = &  \alpha \psi_3(t) - \lambda_3(t),  & & \psi_3(T) = \mu_3,
\end{array} \]
and the switching functions \eqref{switch} must satisfy the conditions:
\[ \begin{array}{ll}
\theta_u (t) = A \psi_2(t) + \psi_3(t) - K = 0, & \theta_w (t) = p - \psi_3(t) > 0, \\
\theta_v (t) = - \psi_2 - 1 \left\{ \begin{array}{rl} < 0, & t \in [0, t_D), \\
= 0, & t \in [t_D,T]. \end{array} \right. &  \end{array} \]
By mean of the same technique employed in previous case-studies, we have found a set of multipliers  $\lambda_i(t) \geq 0$ and $\mu_i \geq  0, I = 1,2,3,4$ that satisfy the complementary slackness condition \eqref{slackness}:
\[ \begin{array}{cccc}
\lambda_1(t) = 0, & \lambda_2(t) = \left\{ \begin{array}{ll}  0, & t \in [0,t_D) \\ r, & t \in [t_D,T] \end{array} \right., & \lambda_3(t) = \alpha (A + K), & \lambda_4(t) = 0 \\
\mu_1 = 0, & \mu_2 = 0, & \mu_3 = A + K, & \mu_4 = 0. \end{array} \]
The latter proves the optimality of control \eqref{oc-scenario3}. In order to find the switching point $t_D \in (0,T),$ we integrate \eqref{system-D} under optimal control \eqref{oc-scenario3} and obtain
\[ D(t) = D_0 \exp[rt] + \frac{1}{r} \left( v_{max} - A w_{max} \right) \left( 1 - \exp[rt] \right), \quad t \in [0, t_D) \]
where $t_D$ indicate the exact moment when $D(t)$ hits zero, that is, $D(t_D)=0.$ Thus,
\begin{equation}
\label{t-d-sc3}
t_D = \frac{1}{r} \ln \frac{v_{max} - A w_{max}}{v_{max} - r D_0 - A w_{max}}.
\end{equation}
Naturally, for smaller initial debt $D_0$ this point will be closer to zero, and for larger one it will be farther from zero.
\begin{remark}
\label{rem3}
Eventually, it may occur that $t_D \geq T$. In this case, the initial debt must be substantially large: $D_0 \geq (v_{max} - A w_{max})(1 - \exp [-rT]).$ Therefore, it will not be totally paid off by final time $T$ even under constant debt repayment at maximum rate:
\[ u^{*} = w_{max}, \quad  v^{*} =  v_{max},  \quad w^{*} = w_{max}. \]
\end{remark}
Clearly, if $t_D < T$, then $D(T)=0$ and the value of criterion \eqref{objective} is solely defined by the cumulative net profit at final time $T$. Direct integration of \eqref{system-N} under optimal control \eqref{oc-scenario3} over $[0,t_D] \cup [t_D,T]$ results in
\begin{equation}
\label{objective-sc3}
 \mathcal{J} ( u^{*},v^{*},w^{*} ) = N_0 +  (A w_{max} - v_{max})t_D + w_{max} (p - A -K) T - BT.
\end{equation}
\begin{remark}
\label{rem4}
Actually, Scenario \ref{subsec33} can be treated as a special case of Scenario \ref{subsec32} when $S_0=0$ and hence $t_S =0 > t_D, S(t) \equiv 0$ for all $t \in [0,T].$ In this case, \eqref{t-d-sc2b} coincides with \eqref{t-d-sc3} and \eqref{objective-sc2} becomes \eqref{objective-sc3}.  On the other hand, Scenario \ref{subsec31} \emph{cannot} be treated as a special case of Scenario \ref{subsec32} by merely setting $t_D=0$. The latter becomes obvious by setting $t_D$ in \eqref{objective-sc2} and comparing this result with \eqref{objective-sc1}.
\end{remark}
In all three scenarios considered above, the optimal control strategies result in total  absence of debts in the end of period, that is $D(T)=0$ (except the situation when initial arrears $D_0$ are extremely high, see Remarks \ref{rem2} and \ref{rem3}). On the other hand, positive cumulative profit and, consequently, positive value of the objective functional $  \mathcal{J} ( u^{*},v^{*},w^{*} )$ can only be guaranteed in Scenario \ref{subsec31} under the ``profitability condition" \eqref{profit-con}. Effectively, second summand in \eqref{objective-sc2} (as well as in \eqref{objective-sc3}) will be negative if $v_{max}$ is very large. The latter may result in negative overall profit $ N(T) <0$ even under the condition \eqref{profit-con}.

Therefore, if there is an initial debt $D_0 > 0$ and available cash $N_0 > 0$ while the firm's capacity of debt repayment is almost unlimited (that is, $v_{max} \to \infty$), one should think about alternative control strategies in order to guarantee the positivity of overall profit.

\section{Alternative control strategies and their optimality}
\label{sec:alter}

Let us consider again the optimal control problem \eqref{objective}-\eqref{cont-constr} under Scenario \ref{subsec32} when the firm has almost unlimited capacity of debt repayment. Mathematically, it means that $v_{max} \to \infty$ and also implies that $t_D \to 0$. In other words, this passage to the limit yields discontinuities in the state variables $N(t)$ and $D(t)$ at the initial point $t=0$:
\[ \lim_{t_D \to 0} D(t_D) = 0 \hspace{1cm} \mbox{(by definition of $t_D$)} \]
\begin{eqnarray*}
\hspace{-10mm} \lim_{t_D \to 0} N(t_D) & = & N_0 + \lim_{\scriptsize \begin{array}{c} v_{max} \to \infty \\ t_D \to 0 \end{array} } \int\limits_{0}^{t_D} \left( p w_{max} - v_{max} - B \right) dt \\
& = & N_0 + \lim_{t_D \to 0} \left( p w_{max} -B \right) t_D - \lim_{\scriptsize \begin{array}{c} v_{max} \to \infty \\ t_D \to 0 \end{array} } v_{max} t_D \\
& = & N_0 - \lim_{v_{max} \to \infty} v_{max} \cdot \frac{1}{r} \ln \frac{v_{max}}{v_{max} - rD_0} = N_0 - \frac{1}{r} \lim_{v_{max} \to \infty} \frac{\ln \frac{\textstyle v_{max}}{\textstyle v_{max} - rD_0}}{\frac{\textstyle 1}{\textstyle v_{max}}} \\
& = & N_0 - \frac{1}{r} \lim_{v_{max} \to \infty} \frac{rD_0 v_{max}^2}{v_{max} (v_{max} - rD_0)} = N_0 - D_0.
\end{eqnarray*}
In the above expressions we have used the form of $t_D$ given by \eqref{t-d-sc2a} (since $t_D < t_S$) together with L'H\^{o}spital's rule. By permitting such finite jumps in the initial states of $N(t)$ and $D(t)$, we can now adjust the optimal control strategy \eqref{oc-scenario2} in a way that its implementation will result in a positive value of the objective functional \eqref{objective}.

From the above formula, it is clear that optimal solution must depend on the relationship between $N_0$ and $D_0$; namely, there are two options to be revised: $N_0 \geq D_0$ and $N_0 < D_0$.

\subsection{Total debt repayment at initial time: $N_0 \geq D_0$}

Suppose that the firm is capable to settle all its debts at once, and still to have a non-negative profit ($N_0 - D_0 \geq 0$) at the beginning of the period. Then it will be profitable first to sell the existing stock of finished goods, without starting the production and, thus, not generating new arrears. Exactly at the moment $t_S$ (when the stock is cleared out, that is, $S(t_S)=0$) the production is started at the rate equal to the maximum volume of demand $w_{max}$. In other words, we arrive to Scenario \ref{subsec31} with one difference only: the initial cash resources are now given as $N_0 - D_0 \geq 0$.  It is easy to prove that control strategy \eqref{oc-scenario1} will be optimal and to do so one should merely repeat all the deductions made in Subsection \ref{subsec31} with $N_0 - D_0$ instead of $N_0$.

Finally, the objective functional will have positive value under optimal control \eqref{oc-scenario1}:
\[ \mathcal{J} \left( u^{*},v^{*},w^{*} \right) = N_0 - D_0 + (p w_{max} - B)T + (A + K)w_{max} (t_S - T) > 0 \]
due to the profitability condition \eqref{profit-con}.

\subsection{Partial debt repayment at initial time: $N_0 < D_0$}

This case looks more challenging than the previous one. The firm does not have enough cash to settle all its debts right away. Therefore, it spends all available cash $N_0$ to repay a part of the initial debt $D_0$. Thus, the firms profit at $t=0$ becomes equal to zero, while its initial debt is reduced to $D_0 - N_0 >0.$  If $S_0 >0$ then there should be no production up to the moment $t_S$ that marks a full clearance of the finished goods stock. From $t_S$ the production at a rate  $w_{max}$ (maximum volume of demand) is started. Meanwhile, all the profit obtained from sales is spent on repayment of previous and new debts right up to the moment $t_D$ at which all the debts are paid off, that is, $D(t_D)=0$. Note that new debts are generated from the commencement of production, that is, for $t \geq t_S$. Thus, for all $t \geq t_D$ the firm's disbursements related to the production process will become equal to $A w_{max}$. In consequence, we propose the following optimal control:
\begin{equation}
\label{oc-42}
\hspace{-10mm} u^{*} = \left\{ \begin{array}{rl} 0, & t \in [0, t_S) \\
\!\!\! w_{max}, & t \in [t_S,T] \end{array} \right. \!\!\!, \;\;  v^{*} = \left\{ \begin{array}{rl} \!\!\! p w_{max} - K u^{*} - B, & t \in [0, t_D) \\
A u^{*}, & t \in [t_D,T] \end{array} \right.\!\!\!, \;\; w^{*} = w_{max}
\end{equation}
that have more sophisticated structure since $v^{*}$ depends also on $u^{*}$.

Lagrange multipliers $\lambda_i(t) \geq 0$ and $\mu_i \geq 0, i=1,2,3,4$ must satisfy the complementary slackness conditions \eqref{slackness}:
\[ N(t) \left\{ \begin{array}{rl} = 0, & t \in [0, t_D) \\ > 0, & t \in [t_D, T] \end{array} \right. \quad \mbox{implies} \quad \lambda_1(t) \left\{ \begin{array}{rl} \geq 0, & t \in [0, t_D) \\ = 0, & t \in [t_D, T] \end{array} \right., \; \mu_1 = 0; \]
\[ D(t) \left\{ \begin{array}{rl} > 0, & t \in [0, t_D) \\ = 0, & t \in [t_D, T] \end{array} \right.  \quad \mbox{implies} \quad \lambda_2(t) \left\{ \begin{array}{rl} = 0, & t \in [0, t_D) \\ \geq 0, & t \in [t_D, T] \end{array} \right., \; \mu_2 \geq 0; \]
\[ S(t) \left\{ \begin{array}{rl} > 0, & t \in [0, t_S) \\ = 0, & t \in [t_S, T] \end{array} \right.  \quad \mbox{implies} \quad \lambda_3(t) \left\{ \begin{array}{rl} = 0, & t \in [0, t_S) \\ \geq 0, & t \in [t_S, T] \end{array} \right., \; \mu_3 \geq 0; \]
\[ S(t) < S_{max} \quad \mbox{implies} \quad \lambda_4 (t) = 0, \; \mu_4 = 0. \]
The underlying adjoint system \eqref{adsystem} can be written as
\[ \begin{array}{rclcl}
\dot{\psi}_1(t) & = & \left\{ \begin{array}{ll} - \lambda_1(t), & t \in [0, t_D) \\ 0, & t \in [t_D,T] \end{array} \right., & & \psi_1(T) = 1, \\
\dot{\psi}_2(t) & = & \left\{ \begin{array}{ll} - r \psi_2(t), & t \in [0, t_D) \\ -r \psi_2(t) - \lambda_2(t), & t \in [t_D,T] \end{array} \right., & & \psi_2(T) = \mu_2 - 1, \\
\dot{\psi}_3(t) & = & \left\{ \begin{array}{ll} \alpha \psi_3(t), & t \in [0, t_S) \\ \alpha \psi_3(t) - \lambda_3(t), & t \in [t_S,T] \end{array} \right., & & \psi_3(T) = \mu_3,
\end{array} \]
and the switching functions \eqref{switch} must satisfy the conditions:
\[ \hspace{-10mm} \theta_u (t) = - K \psi_1(t) + A \psi_2(t) + \psi_3(t)  \left\{ \begin{array}{rl} \!\!\!< 0, & t \in [0, t_S) \\
\!\!\!= 0, & t \in [t_S,T] \end{array} \right. \!\!\!, \;\; \begin{array}{r}
\theta_v (t) = - \psi_1(t) - \psi_2 = 0, \\ \theta_w (t) = p \psi_1 (t) - \psi_3(t) > 0. \end{array} \]
Here (as well as in Scenario \ref{subsec32}) we have two switching points $t_S$ and $t_D$; therefore, the optimal solution will essentially depend on their position with respect to each other. Consequently, one should revise two cases ($t_S < t_D$ and $t_S > t_D$) and find only \emph{one} underlying set of multipliers $\lambda_i(t), \mu_i, i = 1,2,3,4$ that satisfy the complementary slackness conditions \eqref{slackness} in both cases. The latter results in the following set:
\[ \lambda_1(t) = \left\{ \begin{array}{ll} r \exp[r(t_D -t)], & t \in [0, t_D) \\ 0, & t \in [t_D,T] \end{array} \right.\!\!\!, \quad  \lambda_2(t) = \left\{ \begin{array}{ll}  0, & t \in [0,t_D) \\ r, & t \in [t_D,T] \end{array} \right. \!\!\!, \]
\[ \lambda_3(t) = \left\{ \begin{array}{ll}  0, & t \in [0,t_S) \\ \alpha (A + K), & t \in [t_S,T] \end{array} \right.\!\!\!, \quad \begin{array}{l} \lambda_4(t) =0, \\ \mu_1 = \mu_2 = \mu_4 = 0, \;\; \mu_3 = A + K. \end{array} \]
Existence of the above multipliers clearly proves the optimality of control \eqref{oc-42}.
This case bears strong resemblance to the Scenario \ref{subsec32} since there are two switching points ($t_S$ and $t_D$) with the same meanings $S(t_S)=0$ and $ D(t_D)=0$, respectively. Their position with respect to each other will naturally depend on the amount of initial debt $D_0 - N_0 >0$. Eventually, it may also happen that $t_D = t_S$ and there will only one switching point. The following proposition extends the Proposition \ref{prop1} for this case and provides exact formulae for calculation of $t_D$ and $t_S$ in terms of problem entries.
\begin{prop}
\label{prop3}
The stock of finished goods becomes empty at the moment $t_S$ given by \eqref{t-s-sc1} independently of initial debt amount $D_0 - N_0 > 0$. The total debt repayment occurs at the moment $t_D$ that depends on $D_0 - N_0$ in the following way: \\
(a) if $D_0 - N_0 = (p w_{max} - B) \left( 1 - \exp[- rt_S] \right)/r$ then $t_D = t_S$ and there is only one switching point; \\
(b) if $D_0 - N_0 < (p w_{max} - B) \left( 1 - \exp[- rt_S] \right)/r$ then $t_D < t_S$ and
\begin{equation}
\label{t-d-sc4a}
t_D = \frac{1}{r} \ln \left( \frac{p w_{max} - B}{p w_{max} - B - r (D_0 - N_0)} \right);
\end{equation}
(c) if $D_0 - N_0 > (p w_{max} - B) \left( 1 - \exp[- rt_S] \right)/r$ then $t_D > t_S$ and
\begin{equation}
\label{t-d-sc4b}
t_D = \frac{1}{r} \ln \left( \frac{w_{max} (p - A - K) - B}{p w_{max} - B - r (D_0 - N_0) - (A + K) w_{max} \exp[-r t_S]} \right).
\end{equation}
\end{prop}
Formal proof of this proposition can be consulted in the Appendix.
\begin{remark}
\label{rem5}
In effect, if it occurs that $t_S \geq T$ in \eqref{t-s-sc1}, then optimal control \eqref{oc-42} simply becomes
\[ u^{*} = 0, \quad v^{*} = \left\{ \begin{array}{rl} \!\!\! p w_{max} - B, & t \in [0, t_D) \\
0, & t \in [t_D,T] \end{array} \right.\!\!\!, \quad w^{*} = w_{max} \]
and implies no production, only sales of existent stock of finished goods until final time $T$ together with debt repayment up to the moment $t_D$ (if $t_D < T$) or up to $T$ (if $t_D > T$).
\end{remark}
Optimal control strategy \eqref{oc-42} was designed as an alternative to \eqref{oc-scenario2} in order to guarantee the positivity of $ \mathcal{J}^{*} = \mathcal{J} (u^{*},v^{*},w^{*})$. The following proposition provides a sufficient condition under which $ \mathcal{J}^{*} > 0$ regardless of the position of $t_S$ with respect to $t_D$ and $T$.
\begin{prop}
\label{prop4}
If the switching point $t_D$ calculated by either \eqref{t-d-sc4a} or \eqref{t-d-sc4b} satisfies $t_D < T$ then $ \mathcal{J} (u^{*},v^{*},w^{*}) > 0$ where $(u^{*},v^{*},w^{*})$ is given by \eqref{oc-42}. Namely,
\[  \hspace{-10mm} \mathcal{J}^{*} = \left\{ \begin{array}{lcc}
\left[ \left( p - A - K \right) w_{max} - B \right] (T - t_D), & \mbox{if} & 0 < t_S < t_D < T \\
\left[ \left( p - A - K \right) w_{max} - B \right] (T - t_D) + K w_{max} (t_S - t_D), & \mbox{if} & 0 < t_D \leq t_S \leq T \\
\left[ \left( p - A \right) w_{max} - B \right] (T - t_D), & \mbox{if} & 0 < t_D < T < t_S \end{array} \right.\!\!\!. \]
In other words, $\mathcal{J}^{*} >0$ due to the condition \eqref{profit-con}.
\end{prop}
The proof of this proposition is rather straightforward and its key features a given in the Appendix.

Actually, optimal control policies for short-term planning obtained in the Sections \ref{sec:scenarios} and \ref{sec:alter} can be naturally combined by a firm manager or entrepreneur into longer-term decision chains. Effectively, by considering
\[ [0, T] = [T_0, T_1] \cup [T_1, T_2] \cup \cdots \cup [T_{n-1},T_n] \]
with $T_0=0$ and $T_n=T,$ one can successively apply the corresponding scenarios to subintervals $[T_{j-1},T_j], j = 1,2, \ldots,n$ using natural junction conditions for all three state trajectories, that is, terminal conditions at $[T_{j-1},T_j]$ must coincide exactly with initial conditions at $[T_{j},T_{j+1}]$ for all $j= 1,2, \ldots, n-1$.

When ``unconventional'' optimal control policies (introduced in this section) are applied on some $[T_{j-1},T_j]$, the corresponding state trajectories $N^{*}(t)$ and $D^{*}(t)$ will simply have finite jumps at the junction point $T_{j-1}$, while the trajectory $S^{*}(t)$ will remain continuous at $T_{j-1}$. The latter disagrees with the classical theory of optimal control according to which \emph{all} state trajectories must be continuous over the whole interval $[0,T]$. On the other hand, economic meanings of $N(t)$ and $D(t)$ admit such jumps. Therefore, decision chains involving jumps are implementable in practice and they do guarantee (under appropriate initial conditions) a positive outcome as long as the firm remains profitable (that is, while the condition \eqref{profit-con} is kept in force).

\section{Conclusions and future research}
\label{sec:conclusion}

The results of the paper highlight the practical benefits of incorporating ``unconventional'' optimal control strategies in the firm's management in order to corroborate its positive yield and thus to substantiate and stimulate further profitable operation of small business. Such strategies can be easily employed in practice and may provide guidelines to entrepreneurs for short-term and, in form of decision chains, for longer-term planning.

This paper is obviously only a first step toward a more sophisticated modeling and analysis of inventory systems with debt repayment. There are broad perspectives for further research. From the economical point of view, it will be useful to consider this model under market instability with high level of inflation by incorporating an intertemporal discount factor even for short-term planning. In this case, a drastic change of optimal strategies will be quite expectable. It will be also interesting to extend the analysis to systems with more than one product, or with product(s) deterioration and/or with variable rate of demand.

\bibliographystyle{elsarticle-num}
\bibliography{TTV_bib}

\begin{appendix}
\appendix
\renewcommand{\theequation}{A-\arabic{equation}}
  \setcounter{equation}{0}  
  \section*{Appendix: Formal proofs}
\vspace{2mm}

\noindent
\textbf{Proof of the Proposition \ref{prop1}.} The state equation \eqref{system-S} under optimal control strategy \eqref{oc-scenario2} has the same solution as in Scenario \ref{subsec31}; therefore, the switching point $t_S$ is given by \eqref{t-s-sc1}. \\
(a) The situation $t_D = t_S$ implies that all pending debts are repaid just at the same moment as all finished goods are dispatched. On the other hand, optimal strategy \eqref{oc-scenario2} guarantees that $D(t)=S(t)=0$ for all $t \in [t_S,T]$ while
\begin{equation}
\label{D-sc2}
D(t) = \frac{v_{max}}{r} \left( 1 - \exp[rt] \right) + D_0 \exp[rt]
\end{equation}
when $t \in [0,t_S]$. The latter is only possible if $D_0 = v_{max} \left( 1 - \exp[- rt_S] \right)/r$, where $t_S$ is given by \eqref{t-s-sc1}.  \\
(b) For relatively small initial debt (that is, $D_0 < v_{max} \left( 1 - \exp[- rt_S] \right)/r$)
the state trajectory \eqref{D-sc2} hits zero at the moment $t_D$ where $S(t_D) >0, $ that is,  $t_D < t_S$. Therefore, $t_D$ can be calculated as a root of equation $D(t_D)=0:$
\[ t_D = \frac{1}{r} \ln \left( \frac{v_{max}}{v_{max} - r D_0} \right). \]
(c) For rather large initial debt (that is, $D_0 > v_{max} \left( 1 - \exp[- rt_S] \right)/r$) the state trajectory $D(t)$ hits zero at the moment $t_D$ where $S(t_D)=0, $ that is,  $t_D > t_S$. This trajectory has the form of \eqref{D-sc2} for $t \in [0, t_S]$ and vanishes ($D(t)=0$) for $t \in [t_D, T]$. Then, by direct integration of \eqref{system-D} within $[t_S, t_D]$ under initial condition
\[ D(t_S) = D_0 \exp[r t_S] + \frac{v_{max}}{r} \left( 1 - \exp[r t_S] \right) \]
and control strategy \eqref{oc-scenario2}, we obtain that
\[ D(t) =  \left( D_0  + \frac{ A w_{max} \exp[-r t_S] - v_{max}}{r} \right) \exp [rt] - \frac{A w_{max} - v_{max}}{r}, \;\; t \in [t_S, t_D]. \]
Finally, the switching point $t_D$ can be calculated in this case as a root of equation $D(t_D)=0,$ that is,
\[ t_D = \frac{1}{r} \ln \left( \frac{v_{max} - A w_{max}}{v_{max} - r D_0 -A w_{max} \exp[-r t_S]} \right). \]
\vspace{2mm}

\noindent
\textbf{Proof of the Proposition \ref{prop2}.} Effectively, there are two cases: $0 < t_D < t_S < T,$ and $0 < t_S < t_D < T.$

Let $0 < t_D < t_S < T.$ Then by direct integration of \eqref{system-N} under optimal control strategy \eqref{oc-scenario2} with underlying initial conditions it is obtained that
\[ \hspace{-10mm} N^{*}(t) \!= \!\left\{ \!\!\begin{array}{ll}
N_0 + \left( p w_{max} - v_{max} - B \right)t, & t \in [0, t_D] \\
N_0 + \left[ w_{max} (p - A)  - B \right] t +  \left( A w_{max} - v_{max} \right) t_D, & t \in [t_D, t_S] \\
N_0 + \left[ w_{max} (p  - A  - K)  - B \right] t +  \left( A w_{max} - v_{max} \right) t_D + K w_{max} t_S, & t \in [t_S, T]
 \end{array}
\right.\!\!\!. \]
On the other hand, if $0 < t_S < t_D < T$ then direct integration yields:
\[ \hspace{-10mm} N^{*}(t) \!= \!\left\{ \!\!\begin{array}{ll}
N_0 + \left( p w_{max} - v_{max} - B \right)t, & t \in [0, t_S] \\
N_0 + \left( p w_{max} - v_{max} - B \right)t +  K w_{max} (t_S - t), & t \in [t_S, t_D] \\
N_0 + \left[ w_{max} (p - A - K ) - B \right] t +  \left( A w_{max} - v_{max} \right) t_D + K w_{max} t_S, & t \in [t_D, T]
 \end{array}
\right. \!\!\!. \]
In both cases, the value of $\mathcal{J} ( u^{*},v^{*},w^{*} )= N^{*}(T)$ coincides with \eqref{objective-sc2}.
\vspace{2mm}

\noindent
\textbf{Proof of the Proposition \ref{prop3}.} The state equation \eqref{system-S} under optimal control strategy \eqref{oc-42} has the same solution as in Scenarios \ref{subsec31}and \ref{subsec32}; therefore, the switching point $t_S$ is given by \eqref{t-s-sc1}. \\
(a) The situation $t_D = t_S$ implies that all pending debts are repaid just at the same moment as all finished goods are dispatched. On the other hand, optimal strategy \eqref{oc-42} guarantees that $D(t)=S(t)=0$ for all $t \in [t_S,T]$ while
\begin{equation}
\label{D-sc4}
D(t) = (D_0 - N_0) \exp[rt] + \frac{p w_{max} - B}{r} \left( 1 - \exp[rt] \right)
\end{equation}
when $t \in [0,t_S]$. The latter is only possible if $D_0 - N_0 = (p w_{max} - B)\left( 1 - \exp[- rt_S] \right)/r$, where $t_S$ is given by \eqref{t-s-sc1}.  \\
(b) For relatively small initial debt (that is, $D_0 - N_0 < (p w_{max} - B) \left( 1 - \exp[- rt_S] \right)/r)$ the state trajectory \eqref{D-sc4} hits zero at the moment $t_D$ where $S(t_D) >0, $ that is,  $t_D < t_S$. Therefore, $t_D$ is calculated as a root of equation $D(t_D)=0$ with $D(t)$ given by \eqref{D-sc4} what leads to the formula \eqref{t-d-sc4a}. \\
(c) For rather large initial debt (that is, $D_0 - N_0 > (p w_{max} - B) \left( 1 - \exp[- rt_S] \right)/r)$) the state trajectory $D(t)$ hits zero at the moment $t_D$ where $S(t_D)=0, $ that is,  $t_D > t_S$. This trajectory has the form of \eqref{D-sc4} for $t \in [0, t_S]$ and then vanishes ($D(t)=0$) for $t \in [t_D, T]$. Then, by direct integration of \eqref{system-D} within $[t_S, t_D]$ under initial condition
\[ D(t) = (D_0 - N_0) \exp[r t_S] + \frac{p w_{max} - B}{r} \left( 1 - \exp[r t_S] \right) \]
and control strategy \eqref{oc-42}, we obtain that
\begin{eqnarray}
\label{D-t-4}
D(t) & = & (D_0 - N_0) \exp[rt] + \frac{p w_{max} - B}{r} \left( 1 - \exp[rt] \right) \\
& + & \frac{(A + K) w_{max}}{r} \left( \exp[r(t - t_S)] - 1 \right) \nonumber
\end{eqnarray}
for all $ t \in [t_S, t_D].$ Finally, the switching point $t_D$ is calculated in this case as a root of equation $D(t_D)=0$ where $D(t)$ is given by \eqref{D-t-4} what leads to the formula \eqref{t-d-sc4b}.
\vspace{2mm}

\noindent
\textbf{Proof of the Proposition \ref{prop4}.} Effectively, $t_D < T$ implies that all debts have been repaid \emph{before} the end of period. Therefore, $D^{*}(T) = 0$ and $\mathcal{J}^{*} = N^{*}(T)$. The exact form of $N^{*}(T)$ can be obtained by direct integration of \eqref{system-N} under optimal control \eqref{oc-42} using three paths: (a) $0 \rightarrow t_S \rightarrow t_D \rightarrow T; $ (b) $0 \rightarrow t_D \rightarrow t_S \rightarrow T; $ (c) $0 \rightarrow t_D \rightarrow T $ if $t_S > T$.
\end{appendix}

\end{document}